# The central limit theorem for sum-functions of m-tuples of spacings

Sherzod M. Mirakhmedov

V.I. Romanovskiy Institute of Mathematics. Academy of Sciences of Uzbekistan

100174, Tashkent, University str., 9

e-mail: shmirakhmedov@yahoo.com

**Abstract.** Let $n$ points be taken at random on a circle of unit circumference and clockwise ordered. Uniform spacings are defined as the clockwise arc-lengths between the successive points from this sample. We are interested in the behavior of the sums of functions of the $m$-tuples of successive spacings as $n \to \infty$, under the assumption that $m$ can grow together with $n$. We derive asymptotic formulas for the expectation and variance of such sums and establish the asymptotic normality of the sums.

## 1. Introduction

The uniform distribution on the interval [0,1] is of fundamental importance in computer science, probability and statistics. For instance, $F$ is the distribution function of an absolutely continuous random variable (r.v.) $X$ if and only if $F(X)$ has the uniform distribution. Due to this fact, testing the hypothesis that a sample comes from a specified absolutely continuous distribution can be reduced to testing uniformity. One of the most popular approaches to testing uniformity is based on statistics like the sums of spacings, or the gaps between the ordered observations. Besides the hypothesis testing, statistics of this kind are also of great interest in several contexts including reliability (see Pyke (1965)), circular data analysis where they play a pivotal role (see, e.g., Jammalamadaka and SenGupta (2001)), and spacings-based parameter estimation (see, e.g., Ekström et al. (2020). The asymptotic theory for the sums of uniform spacings plays in all these cases an important role (see Pyke (1965), Rao and Sethuraman (1975), Cressie (1976), Kuo and Rao (1981). Hall (1984, 1986), Jammalamadaka et al (1989), L'ecuyer (1997), Mirakhmedov (2010), Mirakhmedov and Jammalamadaka (2013), Ekström et al (2020) . and references therein).

Let $0 \leq U_1 \leq U_2 \leq .... \leq U_{n-1} \leq 1$ be the ordered version of a sample of size $n-1$ from a uniform distribution on [0,1]. We put $U_0 = 0$, $U_n = 1$, and set $U_k = 1 + U_{k-n}$ for $k \geq n$. Two types of spacings are mostly considered in the literature: the overlapping spacings, defined as $S_{k,m} = U_{k+m} - U_k$, and the disjoint spacings, defined by $S_{km,m} = U_{(k+1)m} - U_{km}$, $k = 0,...,n-1$, where $m \geq 1$ is the step, or order, of spacings. The quantities $S_k = S_{k,1}$ are referred to as simple spacings.

The problem we are dealing with in this paper is closely related to the work by Holst (1979), where for a fixed $m \geq 1$ and a given real measurable function $h$ it is proved asymptotic normality of the r.v.'s

$$Z_{n,m} = \sum_{k=0}^{n-1} h(nS_k,...,nS_{k+m-1}). \tag{1.1}$$

This is an extension of the case $m=1$ that had been considered earlier by several authors, see e.g. Le Cam (1958), Rao and Sethuraman (1975) and references therein. The class of statistics (1.1) includes the sums of overlapping spacings, viz.

$$V_{n,m} = \sum_{k=0}^{n-1} h(nS_{k,m}). \tag{1.2}$$

The asymptotic normality of statistics like $V_{n,m}$ has been studied in the past. Beirlant et al. (1991) treated the case of a fixed $m \geq 1$, Creesie (1976) allowed a growing $m = m(n)$ such that $m(n) = o(n^{1/3})$, whereas Hall (1986) assumed that $m(n) = o(n)$ and $n^{1/r} = o(m)$, where $r \geq 2$ is such that the derivative $h^{(r+1)}$ exists and is continuous in a neighborhood of $1$. Creesie (1976) and Hall (1984) considered also the sums of spacings (based on the interval $]0,1]$) of the form

$$W_{n,m} = \sum_{k=0}^{n-m-1} h(nS_{k,m}). \tag{1.3}$$

However, Hall (1986) and Guttorp and Lockhart (1989) pointed out some shortcomings of this statistic and recommend using statistics on a circle rather than on the line. We note that, under a mild condition, the r.v.'s $V_{n,m}$ and $W_{n,m}$ have the same asymptotically normal distribution, see Corollary 2.3 below.

The statistics $Z_{n,m}$, $V_{n,m}$ and $W_{n,m}$ are symmetric in spacings. In this paper, we consider a general class of statistics of the form

$$R_{n,m} = \sum_{k=0}^{n-1} h_k(nS_k,...,nS_{k+m-1}), \tag{1.4}$$

where $h_k$s are a given sequence of real measurable functions. In particular, if $h_k = 0$ for $k = n-m+1,...,n-1$ we obtain the statistic

$$R^*_{n,m} = \sum_{k=0}^{n-m} h_k(nS_k,...,nS_{k+m-1}). \tag{1.5}$$

Clearly, the r.v.'s $Z_{n,m}$, $V_{n,m}$ and $W_{n,m}$ are special cases of $R_{n,m}$. By setting $h_k(x_1,...,x_m) = h_k(x_1 + ... + x_m)$ for $k = lm$, where $l = 0,1,...,[n/m]-1$, and otherwise $h_k(x_1,...,x_m) = 0$, we get the sums of disjoint spacings, viz.,

$$Q_{n,m} = \sum_{k=0}^{[n/m]} h_k(nS_{km,m}). \tag{1.6}$$

This statistic was suggested by Del Pino (1979) who examined its behavior for a fixed *m*. For a detailed asymptotic theory for r.v.'s (1.6), including the case where *m* increases as the sample size grows, we refer to Jammalamadaka et al (1989), Guttorp and Lockhart (1989), Mirakhmedov (2005, 2010), Mirakhmedov et al (2011), Mirakhmedov and Jammalamadaka2013) and Ekström et al (2020) .

We would like to stress that non-symmetric statistics are useful, among other applications, when studying the efficiencies of spacings-based tests for uniformity. Namely, statistics based on the spacings generated by sequences of local alternative distributions converging to the uniform one can be reduced to a non-symmetric version of uniform spacings-statistics, see, for instance, Kio and Jammalamadaka (1981). We remind that, say, in the context of spacings-based parameter estimation as well as in the problem of testing for uniformity, if we are dealing with spacings of the same order then statistics based on overlapping spacings are considered preferable to statistics based on disjoint spacings. However, working with statistics like $Q_{n,m}$ involves considerably less calculation. Moreover, in testing for uniformity versus to the sequences of close alternatives, by choosing the value of spacings order *m* corresponding to $Q_{n,m}$ larger than the one corresponding to $V_{n,m}$, the $Q_{n,m}$-test can be made as efficient as $V_{n,m}$-test, see Jammalamadaka et al (1989). The asymptotic theory here suggests that larger *m* values are always better.

It is the purpose of this paper to establish asymptotically normality of the r.v. $R_{n,m}$ and provide asymptotic expressions for its expectation and variance. We include situations where $m = m(n) \to \infty$ as $n \to \infty$, with the proviso that $m = o(n)$. Reviewing the available literature it was noticed that Holst's formula for asymptotic variance needs to be adjusted. As consequence we provide corrected and extended for increasing *m* version of Holst (1979) theorem. Next, we apply general results to the $V_{n,m}$ statistic and its three important variants; the consequences, in turn, are a complement and refinement of previously known results.

The spacings clearly depend on *n*, and the functions *h* and $h_k$ may also depend on *n* in the cases we consider, but we will suppress the extra subscript *n* to make our notation less cumbersome. Throughout the paper, *c* denotes a positive absolute constant, possibly not the same one in different occurrences, every unspecified asymptotic relation considered is as $n \to \infty$, $d_n \sim b_n$ stands for the $d_n / b_n \to 1$, $\xrightarrow{d}$ denotes convergence in distribution, $\mathcal{N}(0,1)$ stands for the standard normal r.v., for given quantities $a_i$ we set $a_k^m := (a_k,...,a_{k+m-1})$, $na_k^m := (na_k,...,na_{k+m-1})$, $|a_k^m| := a_k + ... + a_{k+m-1}$. Section 2 states our main results, their proofs are presented in Section 3.

## 2. Results.

Let $X_0, X_1, ..., X_{n-1}$ be independent r.v.s having a common standard exponential distribution. We set $X_{n+j} = X_j$, $\bar{h}_k(X_k^m) = \bar{h}_k(X_k^m) - Eh_k(X_k^m)$.

**Proposition 2.1.** Let $E|h_k(X_k^m)(|X_k^m| - m)^3| < \infty$, $k = 0, 1, ..., n-1$. Then

$$ER_{n,m} = A_{n,m} + \frac{1}{2n} \sum_{k=0}^{n-1} \text{cov}\left(h_k(X_k^m), (|X_k^m| - m) - (|X_k^m| - m)^2\right)(1 + o(1)) + r_{1n} \ ,$$

where

$$A_{n,m} = \sum_{k=0}^{n-1} Eh_k(X_k^m) ,$$

$$|r_{1n}| \leq \frac{c}{n^{3/2}} \sum_{k=0}^{n-1} E|\bar{h}_k(X_k^m)(|X_k^m| - m)^3| \leq \frac{c_1 m^{3/2}}{n^{3/2}} \sum_{k=0}^{n-1} \sqrt{\text{var} \, h_k(X_k^m)} \ ,$$

provided $Eh_k^2(X_k^m) < \infty$ for the last inequality.

**Proposition 2.2.** Assume $E|h_k(X_k^m)|^3 < \infty$, $k = 0, 1, ..., n-1$. Then

$$\text{Var} R_{n,m} = \sigma_{n,m}^2 (1 + r_{2n})$$

where $\sigma_{n,m}^2 = n(C_{n,m} - B_{n,m}^2)$,

$$C_{n,m} = \frac{1}{n} \sum_{k=0}^{n-1} \sum_{j=k-m+1}^{k+m-1} \text{cov}\left(h_k(X_k^m), h_j(X_j^m)\right), \quad B_{n,m} = \frac{1}{n} \sum_{k=0}^{n-1} \text{cov}\left(h_k(X_k^m), |X_k^m|\right) , \quad (2.1)$$

$$|r_{2n}| \leq \frac{c}{\sigma_{n,m}^2} \left[ \frac{m}{n} \sum_{k=0}^{n-1} \sum_{j=k-m+1}^{k+m-1} \left(E|\bar{h}_k(X_k^m)|^3 E|\bar{h}_j(X_j^m)|^3\right)^{1/3} + \frac{m^{3/2}}{n^{3/2}} \left(\sum_{k=0}^{n-1} (E|\bar{h}_k(X_k^m)|^3)^{1/3}\right)^2 \right] .$$

Set $g_{k,n}(X_k^m) = h_k(X_k^m) - Eh_k(X_k^m) - B_{n,m}(X_k - 1)$.

**Theorem 2.1.** Let $m = m(n) = o(n)$ and the sequence of functions $h_0, ..., h_{n-1}$ is such that

(i) $C_{n-m,m} \sim C_{n,m}$ and $B_{n-m,m} \sim B_{n,m}$,

(ii) for some $r > 2$

$$\frac{m^{r-1}}{\sigma_{n,m}^r} \sum_{k=0}^{n-1} E|g_{k,n}(X_k^m)|^r \to 0 . \quad (2.2)$$

Then $(R_{n,m} - A_{n,m})/\sigma_{n,m} \xrightarrow{d} \mathcal{N}(0,1)$ .

As immediate consequence of this theorem we have

**Theorem 2.2.** If $m = m(n) = o(n)$, and for some $r > 2$

$$\frac{m^{r-1}}{\sigma_{n-m,m}^r} \sum_{k=0}^{n-m-1} E|g_{k,n-m}(X_k^m)|^r \to 0 ,$$

then $(R_{n,m}^* - A_{n-m,m})/\sigma_{n-m,m} \xrightarrow{d} \mathcal{N}(0,1)$ .

Let's consider the "symmetric" statistics $Z_{n,m}$ (1.1). We have

$$A_{n,m} = nEh(X_0^m), \quad \sigma_{n,m}^2 = n\left[\sum_{j=-m+1}^{m-1} \text{cov}\left(h(X_0^m), h(X_j^m)\right) - \text{cov}^2\left(h(X_0^m), |X_0^m|\right)\right] =: n\sigma_m^2. \qquad (2.3)$$

**Corollary 2.1.** Let $Eh^2(X_0^m) < \infty$, then

$$EZ_{n,m} = nEh(X_0^m) + \frac{1}{2}\text{cov}\left(h(X_0^m), (|X_0^m| - m) - (|X_0^m| - m)^2\right)(1+o(1)) + \frac{cm^{3/2}}{\sqrt{n}}\sqrt{\text{var}\, h(X_0^m)}$$

**Corollary 2.2.** If $E|h(X_0^m)|^3 < \infty$ then

$$\text{Var}Z_{n,m} = n\sigma_m^2\left(1 + c\frac{m^{3/2}\left(E|\bar{h}(X_0^m)|^3\right)^{2/3}}{\sqrt{n}\sigma_m^2}\right).$$

Set $g(X_0^m) = h(X_0^m) - Eh(X_0^m) - (X_0 - 1)\text{cov}\left(h(X_0^m), |X_0^m|\right)$.

**Corollary 2.3.** Let $m = m(n) = o(n)$ and for some $r > 2$

$$\frac{m^{r-1}}{n^{(r-2)/2}\sigma_m^r} E|g(X_0^m)|^r \to 0.$$

Then $(Z_{n,m} - nEh(X_0^m))/\sigma_m\sqrt{n} \xrightarrow{d} \mathcal{N}(0,1)$.

The same assertion hold for the r.v.

$$Z_{n,m}^* = \sum_{k=0}^{n-m-1} h(nS_k, \ldots, nS_{k+m-1}).$$

That is the r.v.s $Z_{n,m}$ and $Z_{n,m}^*$ have the same asymptotically normal distribution.

**Remark 2.1.** Holst (1979) stated that if $Eh(X_0^m) = 0$, $0 < \text{Var}\,h(X_0^m) < \infty$ and $m$ is fixed then $n^{-1/2}Z_{n,m} \xrightarrow{d} \mathcal{N}(0, \sigma^2)$, where

$$\sigma^2 := \sum_{j=-m+1}^{m-1} \text{cov}\left(h(X_0^m), h(X_j^m)\right) - \left(m^{-1}\sum_{j=-m+1}^{m-1} \text{cov}\left(h(X_0^m), |X_j^m|\right)\right)^2.$$

Corollary 2.2 shows that this formula for $\sigma^2$ is not correct and one should use the adjusted variant given by $\sigma_m^2$ in (2.3). What happened was that, in the course of proof in Holst (1979), it was assumed that two random variables considered on p.1068 (lines 8 and 10) have the same asymptotic distribution. However, they have asymptotically normal distributions with different variances. Corrections must be made for the variance of the sum of spacings used in other papers as well, for example, the formula for $\sigma_m^2$ in L'ecuyer (1997), p. 254.

Further, we consider applications of Corollary 2.3 to symmetric statistic $V_{n,m}$ (1.2) and its three important versions. The r.v. $V_{n,m}$ is special case of $R_{n,m}$ where $h_k(x_k^m) = h(|x_k^m|)$, $h(u)$ is a real measurable function defined on non-negative axis, therefore as a straightforward consequence of Corollary 2.3 we have

**Corollary 2.4.** Let $m = m(n) = o(n)$ and for some $r > 2$

$$\frac{m^{r-1}}{n^{(r-2)/2}\sigma_m^r} E|g(|X_0^m|)|^r \to 0.$$

Then $(V_{n,m} - nEh(|X_0^m|))/\sigma_m\sqrt{n} \xrightarrow{d} \mathcal{N}(0,1)$, where

$$g(|X_0^m|) = h(|X_0^m|) - Eh(|X_0^m|) - (X_0 - 1)\text{cov}\left(h(|X_0^m|), |X_0^m|\right),$$

$$\sigma_m^2 = \sum_{j=-m+1}^{m-1} \text{cov}\left(h(|X_0^m|), h(|X_j^m|)\right) - \text{cov}^2\left(h(|X_0^m|), |X_0^m|\right). \tag{2.4}$$

This Corollary 2.4 extends Theorem of Beirlant et al (1991) for the case $m = m(n) \to \infty$. It is also complement to Theorem 2 of Hall (1986) covering the case where $m(n) = O(n^{1/r})$, see above Introduction. Most common examples of $V_{n,m}$ are the so-called Greenwood statistic

$$G_{n,m}^2 = \sum_{k=0}^{n-1}(nS_{k,m})^2,$$

the Log-spacings statistic, which is sometimes called the Moran statistic,

$$M_{n,m} = \sum_{k=0}^{n-1}\ln(nS_{k,m}),$$

the entropy-type statistic

$$H_{n,m} = \sum_{k=0}^{n-1}(nS_{k,m})\ln(nS_{k,m}).$$

Let $\psi(m)$ and $\zeta(k,m)$ stands for digamma and Hurwitz zeta function, respectively, viz.,

$$\psi(m) = \sum_{k=1}^{m-1}\frac{1}{k} - \gamma, \quad \zeta(k,m) = \sum_{j=m}^{\infty} j^{-k}, \quad \gamma = 0.5772... \text{ is Euler's constant.}$$

**Corollary 2.5 .**

(i) If $m = o(n)$ then $\sqrt{3}\left(G_{n,m} - nm(m+1)\right)/\sqrt{2m(m+1)(2m+1)n} \xrightarrow{d} \mathcal{N}(0,1)$, (2.5)

(ii) If $m = o(n^{1/3})$ then

$$\left(M_{n,m} - n\psi(m)\right)/\sqrt{n\left((2m^2 - 2m + 1)\zeta(2,m) - 2m + 1\right)} \xrightarrow{d} \mathcal{N}(0,1), \tag{2.6}$$

and

$$\left(H_{n,m} - nm\psi(m+1)\right)/\sqrt{n\left(2(m(m+1))^2\zeta(2,m+2) - m(m+1)(2m-1)\right)/4} \xrightarrow{d} \mathcal{N}(0,1). \tag{2.7}$$

**Remark 2.2**. Condition $m = o(n^{1/3})$ is due to fact that we applied Corollary 2.4 with $r = 4/$ This is because we get an upper bound for $Eg^4(|X_0^m|)$ only, which already required too lengthy calculations.

We conjecture that asymptotic normality of $M_{n,m}$ and $H_{n,m}$ hold for $m=o(n^{(r-2)/2(r-1)})$, arbitrary integer $r \geq 3$, the condition close to $m=o(n^{1/2})$.

**Remark 2.3**. On using Theorem 2.2 (by setting $h_k(x_k^m)=h_k(|x_k^m|), k=0,1,...,n-1$) and Corollary 2.5 one can extend results of section 5 of Misra and van der Meulen (2001) for the case where the order of spacins $m \to \infty$, $m=o(\sqrt{n})$. In particular, for this case assertion (2.7) is extension of their Corollary 4.2.

## 3. Proofs

The following lemma is a straightforward consequence of the formula (1.2) by Mirakhmedov (2005). We will use it in various forms.

**Lemma**. Assume the function $T(x_0,...,x_{n-1})$ is such that

$$\int_{-\infty}^{\infty} \left| E\left(T(X_0^n) \exp\{i\tau |X_0^n|\}\right) \right| d\tau < \infty.$$

Then

$$ET(nS_0^n) = \frac{1}{\mu_n} \int_{-\infty}^{\infty} E\left(T(X_0^n) \exp\{i\tau n^{-1/2}(|X_0^n|-n)\}\right) d\tau,$$

where

$$\mu_n = 2\pi \sqrt{n} n^{n-1} e^{-n} / (n-1)! = \sqrt{2\pi}(1+o(1)). \tag{3.1}$$

Let $\langle w \rangle$ stands for the number of elements of the set $w$, $J_{k,m}=(k,...,k+m-1)$ and

$$w_{k,j} = J_{k,m} \cup J_{j,m}, \quad \bar{w}_{k,j} = (0,1,...,n-1) - w_{k,j}, \quad \phi_{n,m}^{(k,j)}(\tau) = E \exp\left\{i\tau n^{-1/2} \sum_{v \in \bar{w}_{k,j}} (X_v - 1)\right\}.$$

We have

$$\left|\phi_{n,m}^{(k,j)}(\tau)\right| = \left(1+\frac{\tau^2}{n}\right)^{(n-\langle w_{k,j}\rangle)/2}. \tag{3.2}$$

Also Lemma 4 of Petrov (1975, p.140) and that $\langle w_{k,j} \rangle \leq 2m-1 = o(n)$ yields for $|\tau| \leq \sqrt{n}/3$

$$\left|\phi_{n,m}^{(k,j)}(\tau) - \exp\left\{-\frac{\tau^2(n-\langle w_{k,j}\rangle)}{2n}\right\}\left(1+\frac{(i\tau)^3(1-n^{-1}\langle w_{k,j}\rangle)}{6\sqrt{n}}\right)\right| \leq \frac{c\tau^4}{n}\exp\left\{-\frac{\tau^2}{12}\right\}. \tag{3.3}$$

**Proof of Proposition 2.1**. Recall notation $\bar{h}_k(x_k^m) = h_k(x_k^m) - Eh_k(X_k^m)$. We have

$$\int_{-\infty}^{\infty} \left|E\left[\bar{h}_k(X_k^m)\exp\{i\tau(|X_0^n|-n)\}\right]\right| d\tau = \int_{-\infty}^{\infty}\left|E\left[\bar{h}_k(X_k^m)\exp\{i\tau(|X_k^m|-m)\}\right]\phi_{n,m}^{(k,k)}(\tau\sqrt{n})\right| d\tau$$

$$\leq E|\bar{h}_k(X_k^m)(X_k^m)| \int_{-\infty}^{\infty}\left|\phi_{n,m}^{(k,k)}(\tau\sqrt{n})\right| d\tau < \infty.$$

since (3.2), $\langle w_{k,k} \rangle = m$ and $m=o(n)$. By Lemma we have the following chine of equalities.

$$ER_{n,m} - A_{n,m} = \mu_n^{-1} \sum_{k=0}^{n-1} \int_{-\infty}^{\infty} E\bar{h}_k(X_k^m) \exp\left\{i\tau n^{-1/2}(|X_0^n| - n)\right\} d\tau$$

$$= \mu_n^{-1} \sum_{k=0}^{n-1} \int_{-\infty}^{\infty} E\left[\bar{h}_k(X_k^m) \exp\left\{i\tau n^{-1/2}(|X_k^m| - m)\right\}\right] \phi_{n,m}^{(k,k)}(\tau) d\tau$$

$$= \frac{1}{\mu_n \sqrt{n}} \sum_{k=0}^{n-1} E\bar{h}_k(X_k^m)(|X_k^m| - m) \int_{-\infty}^{\infty} (i\tau) \exp\left\{-\frac{\tau^2(n-m)}{2n}\right\}\left(1 + \frac{(i\tau)^3(1 - mn^{-1})}{6\sqrt{n}}\right) d\tau$$

$$- \frac{1}{2n\mu_n} \sum_{k=0}^{n-1} E\left[\bar{h}_k(X_k^m)(|X_k^m| - m)^2\right] \int_{-\infty}^{\infty} \tau^2 \exp\left\{-\frac{\tau^2(n-m)}{2n}\right\}\left(1 + \frac{(i\tau)^3(1 - mn^{-1})}{6\sqrt{n}}\right) d\tau$$

$$+ \mu_n^{-1} \sum_{k=0}^{n-1} \int_{-\infty}^{\infty} E\left[\bar{h}_k(X_k^m)\left[\exp\left\{\frac{i\tau}{\sqrt{n}}(|X_k^m| - m)\right\} - 1 - \frac{i\tau}{\sqrt{n}}(|X_k^m| - m) - \frac{(i\tau)^2}{2n}(|X_k^m| - m)^2\right]\right]$$

$$\cdot \exp\left\{-\frac{\tau^2(n-m)}{2n}\right\}\left(1 + \frac{(i\tau)^3(1 - mn^{-1})}{6\sqrt{n}}\right) d\tau$$

$$+ \mu_n^{-1} \sum_{k=0}^{n-1} \int_{-\infty}^{\infty} E\left[\bar{h}_k(X_k^m)\left[\exp\left\{\frac{i\tau}{\sqrt{n}}(|X_k^m| - m)\right\} - 1\right]\right]$$

$$\cdot \left[\phi_{n,m}^{(k,k)}(\tau) - \exp\left\{-\frac{\tau^2(n-m)}{2n}\right\}\left(1 + \frac{(i\tau)^3(1 - mn^{-1})}{6\sqrt{n}}\right)\right] d\tau.$$

Now using (3.3) by quite obvious calculations taking into account (3.1) we obtain

$$ER_{n,m} = A_{n,m} + \frac{1}{2n} \sum_{k=0}^{n-1} E\left[\bar{h}_k(X_k^m)\left((|X^m| - m) - (|X^m| - m)^2\right)\right](1 + o(1))$$

$$+ \frac{c}{n^{3/2}} \sum_{k=0}^{n-1} E\left|\bar{h}_k(X_k^m)\right|\left|X_k^m| - m\right|^3 + \left\||X_k^m| - m\right\|.$$

The second inequality for $|r_{1n}|$ follows by Holder's inequality and that $E(|X_k^m| - m)^6 \leq cm^3$.

**Proof of Proposition 2.2.** By Lemma

$$E(R_{n,m} - A_{n,m})^2 = \mu_n^{-1} \int_{-\infty}^{\infty} E\left[\left(\sum_{k=0}^{n-1} \bar{h}_k(X_k^m)\right)^2 \exp\left\{i\tau n^{-1/2} \sum_{v=0}^{n-1}(X_v - 1)\right\}\right] d\tau =: \mu_n^{-1} \int_{-\infty}^{\infty} I(\tau) d\tau. \quad (3.4)$$

To keep future notation simple we set $Y_{k,j} := \sum_{l \in w_{k,j}} (X_l - 1)$,

$$\alpha_{k,j}(\tau) := E\left[\bar{h}_k(X_k^m)\bar{h}_j(X_j^m) \exp\left\{i\tau n^{-1/2} Y_{k,j}\right\}\right],$$

$$\beta_{k,j}(\tau) := E\bar{h}_k(X_k^m) \exp\left\{i\tau n^{-1/2}(|X_k^m| - m)\right\} E\bar{h}_j(X_j^m) \exp\left\{i\tau n^{-1/2}(|X_j^m| - m)\right\}.$$

Rewrite the integral of (3.4) as follows

$$\int_{-\infty}^{\infty} I(\tau) d\tau = \sum_{k=0}^{n-1} \sum_{j=0}^{n-1} \int_{-\infty}^{\infty} E\left[\bar{h}_k(X_k^m)\bar{h}_j(X_j^m) \exp\left\{i\tau n^{-1/2} \sum_{v=0}^{n-1}(X_v - 1)\right\}\right] d\tau$$

$$= \sum_{k=0}^{n-1} \sum_{j}^{*} \int_{-\infty}^{\infty} \alpha_{k,j}(\tau) \phi_{n,m}^{(k,j)}(\tau) d\tau + \sum_{k=0}^{n-1} \sum_{j}^{**} \int_{-\infty}^{\infty} \beta_{k,j}(\tau) \phi_{n,m}^{(k,j)}(\tau) d\tau =: I_1 + I_2, \quad (3.5)$$

where $\Sigma_j^*$ and $\Sigma_j^{**}$ stand for the summations over $j = 0, \ldots, n-1$ such that $J_{j,m} \cap J_{k,m} \neq \phi$ and $J_{j,m} \cap J_{k,m} = \phi$, respectively. We have

$$I_1 = \sum_{k=0}^{n-1} \sum_{j=k-m+1}^{k+m-1} \int_{-\infty}^{\infty} \alpha_{k,j}(\tau) \exp\left\{-\frac{\tau^2}{2} \frac{n - \langle w_{k,j} \rangle}{n}\right\} \left(1 + \frac{(i\tau)^3 (1 - \langle w_{k,j} \rangle n^{-1})}{6\sqrt{n}}\right) d\tau$$

$$+ \sum_{k=0}^{n-1} \sum_{j=k-m+1}^{k+m-1} \int_{-\infty}^{\infty} \alpha_{k,j}(\tau) \left(\phi_{n,m}^{(k,j)}(\tau) - \exp\left\{-\frac{\tau^2}{2} \frac{n - \langle w_{k,j} \rangle}{n}\right\} \left(1 + \frac{(i\tau)^3 (1 - \langle w_{k,j} \rangle n^{-1})}{6\sqrt{n}}\right)\right) d\tau =: I_{11} + I_{12}. \quad (3.6)$$

We have

$$I_{11} = \sum_{k=0}^{n-1} \sum_{j=k-m+1}^{k+m-1} E\left[\bar{h}_k(X_k^m) \bar{h}_j(X_j^m)\right] \int_{-\infty}^{\infty} \exp\left\{-\frac{\tau^2}{2} \frac{n - \langle w_{k,j} \rangle}{n}\right\} \left(1 + \frac{(i\tau)^3 (1 - \langle w_{k,j} \rangle n^{-1})}{6\sqrt{n}}\right) d\tau$$

$$+ \sum_{k=0}^{n-1} \sum_{j=k-m+1}^{k+m-1} \int_{-\infty}^{\infty} E\left[\bar{h}_k(X_k^m) \bar{h}_j(X_j^m) \left(\exp\{i\tau n^{-1/2} Y_{k,j}\} - 1 - \frac{i\tau}{\sqrt{n}} Y_{k,j}\right)\right]$$

$$\cdot \exp\left\{-\frac{\tau^2}{2} \frac{n - \langle w_{k,j} \rangle}{n}\right\} \left(1 + \frac{(i\tau)^3 (1 - \langle w_{k,j} \rangle n^{-1})}{6\sqrt{n}}\right) d\tau$$

$$+ \frac{1}{6n} \sum_{k=0}^{n-1} \sum_{j=k-m+1}^{k+m-1} E\left[\bar{h}_k(X_k^m) \bar{h}_j(X_j^m) Y_{k,j}\right] (1 - \langle w_{k,j} \rangle n^{-1}) \int_{-\infty}^{\infty} \tau^4 \exp\left\{-\frac{\tau^2}{2} \frac{n - \langle w_{k,j} \rangle}{n}\right\} d\tau.$$

Use the fact that $\langle w_{k,j} \rangle \leq 2m - 1 = o(n)$ and inequality $|e^{it} - 1 - it| \leq t^2/2$ to get

$$\int_{-\infty}^{\infty} I_{11}(\tau) d\tau = \sqrt{2\pi}(1 + o(1)) \sum_{k=0}^{n-1} \sum_{j=k-m+1}^{k+m-1} E\left[\bar{h}_k(X_k^m) \bar{h}_j(X_j^m)\right]$$

$$+ \frac{c}{n} \sum_{k=0}^{n-1} \sum_{j=k-m+1}^{k+m-1} E\left[|\bar{h}_k(X_k^m) \bar{h}_j(X_j^m)| (Y_{k,j}^2 + |Y_{k,j)}|)\right]$$

$$= \sqrt{2\pi}(1 + o(1)) \sum_{k=0}^{n-1} \sum_{j=k-m+1}^{k+m-1} E\left[\bar{h}_k(X_k^m) \bar{h}_j(X_j^m)\right]$$

$$+ \frac{cm}{n} \sum_{k=0}^{n-1} \sum_{j=k-m+1}^{k+m-1} \left(E|\bar{h}_k(X_k^m)|^3 E|\bar{h}_j(X_j^m)|^3\right)^{1/3}, \quad (3.7)$$

here a generalized Hölder's inequality and the facts that $E|Y_{k,j}|^3 \leq (2m)^{3/2}$, $EY_{k,j}^6 \leq cm^3$ are used. In order to estimate integral $I_{12}$ we first split it into two integrals: over $|\tau| \leq \sqrt{n}/3$ and $|\tau| > \sqrt{n}/3$, and next we apply (3.2) and (3.3), respectively. Then by a simple algebra we obtain

$$|I_{12}| \leq \frac{c}{n} \sum_{k=0}^{n-1} \sum_{j=k-m+1}^{k+m-1} E|\bar{h}_k(X_k^m) \bar{h}_j(X_j^m)|$$

$$\leq \frac{c}{n} \sum_{k=1}^{n} \sum_{j=k-m+1}^{k+m-1} \left( E|\bar{h}_k(X_k^m)|^3 E|\bar{h}_j(X_j^m)|^3 \right)^{1/3}, \quad (3.8)$$

by Hölder's inequality and inequalities between moments of a r.v..

Let's consider $I_2$. Note that $\langle w_{k,j} \rangle = 2m$ in the case $J_{j,m} \cap J_{k,m} = \phi$. We have

$$I_2 = \sum_{k=0}^{n-1} \sum_{j}^{**} \int_{-\infty}^{\infty} \beta_{k,j}(\tau) \exp\left\{-\frac{\tau^2}{2}\left(1-\frac{2m}{n}\right)\right\} d\tau$$

$$+ \sum_{k=0}^{n-1} \sum_{j}^{**} \int_{-\infty}^{\infty} \beta_{k,j}(\tau) \left( \phi_{n,m}^{(k,j)}(\tau) - \exp\left\{-\frac{\tau^2}{2}\left(1-\frac{2m}{n}\right)\right\} \right) d\tau = I_{21} + I_{22}. \quad (3.9)$$

Write

$$\sum_{k=0}^{n-1} \sum_{j}^{**} \beta_{k,j}(\tau) = \sum_{k=0}^{n-1} E\bar{h}_k(X_k^m) \exp\left\{\frac{i\tau}{\sqrt{n}}(|X_k^m|-m)\right\} \sum_{j=0}^{n-1} E\bar{h}_j(X_j^m) \exp\left\{\frac{i\tau}{\sqrt{n}}(|X_j^m|-m)\right\}$$

$$- \sum_{k=0}^{n-1} \sum_{j=k-m+1}^{k+m-1} E\bar{h}_k(X_k^m) \exp\left\{\frac{i\tau}{\sqrt{n}}(|X_k^m|-m)\right\} E\bar{h}_j(X_j^m) \exp\left\{\frac{i\tau}{\sqrt{n}}(|X_j^m|-m)\right\} =: V_1(\tau) - V_2(\tau). \quad (3.10)$$

Further

$$V_1(\tau) = \left[ \sum_{k=0}^{n-1} E\bar{h}_k(X_k^m) \left( \exp\left\{\frac{i\tau}{\sqrt{n}}(|X_k^m|-m)\right\} - 1 - \frac{i\tau}{\sqrt{n}}(|X_k^m|-m) \right) + \frac{i\tau}{\sqrt{n}} \sum_{k=0}^{n-1} E\bar{h}_k(X_k^m)(|X_k^m|-m) \right]^2.$$

By this and inequality $|e^{it} - 1 - it| \leq t^2/2$ we obtain

$$\int_{-\infty}^{\infty} V_1(\tau) \exp\left\{-\frac{\tau^2}{2}\left(1-\frac{2m}{n}\right)\right\} d\tau = -\frac{1}{n}\left(\sum_{k=0}^{n-1} E\bar{h}_k(X_k^m)(|X_k^m|-m)\right)^2 (\sqrt{2\pi} + o(1))$$

$$+ \frac{c}{n^2}\left(\sum_{k=0}^{n-1} E|\bar{h}_k(X_k^m)|(|X_k^m|-m)^2\right)^2 + \frac{c}{n^{3/2}} \sum_{k=0}^{n-1} E|\bar{h}_k(X_k^m)|(|X_k^m|-m)^2 \left|\sum_{k=0}^{n-1} E\bar{h}_k(X_k^m)(|X_k^m|-m)\right|$$

$$= -nB_{n,m}^2 \sqrt{2\pi}(1+o(1)) + \frac{cm^{3/2}}{n^{3/2}}\left(\sum_{k=0}^{n-1}(E|\bar{h}_k(X_k^m)|^3)^{1/3}\right)^2. \quad (3.11)$$

since $E|\bar{h}_k(X_k^m)(|X_k^m|-m)^v| \leq c\sqrt{m^v \operatorname{var} h_k(X_k^m)} \leq cm^{v/2}\left(E|\bar{h}_k(X_k^m)|^3\right)^{1/3}$ by Holder's inequality.

Recalling that $E\bar{h}_k(X_k^m) = 0$ it is easy to see that

$$\left|\int_{-\infty}^{\infty} V_2(\tau) \exp\left\{-\frac{\tau^2}{2}\left(1-\frac{2m}{n}\right)\right\} d\tau\right| \leq \frac{c}{n} \sum_{k=0}^{n-1} \sum_{j=k-m+1}^{k+m-1} E|\bar{h}_k(X_k^m)(|X_k^m|-m)| E|\bar{h}_j(X_j^m)(|X_j^m|-m)|$$

$$\leq \frac{cm}{n} \sum_{k=0}^{n-1} \sum_{j=k-m+1}^{k+m-1} \left( E|\bar{h}(X_k^m)|^3 E|\bar{h}_j(X_j^m)|^3 \right)^{1/3}. \quad (3.12)$$

Thus (3.10), (3.11) and (3.12) gives

$$\int_{-\infty}^{\infty} I_{21}(\tau) d\tau = -\sqrt{2\pi}(1+o(1))nB_n^2 + \frac{cm^{3/2}}{n^{3/2}}\left(\sum_{k=0}^{n-1}(E|\bar{h}_k(X_k^m)|^3)^{1/3}\right)^2$$

$$+\frac{cm}{n}\sum_{k=0}^{n-1}\sum_{j=k-m+1}^{k+m-1}\left(E|\bar{h}(X_k^m)|^3 E|\bar{h}_j(X_j^m)|^3\right)^{1/3}. \tag{3.13}$$

By similar reasoning as above using $E\bar{h}_l(X_l^m)=0$, (3.2) and (3.3) one can observe that

$$\left|\int_{-\infty}^{\infty} I_{22}(\tau)d\tau\right| \leq \frac{c}{n^{3/2}}\left(\sum_{k=0}^{n-1} E|\bar{h}_k(X_k^m)(|X_k^m|-m)|\right)^2$$

$$+\frac{cm}{n}\sum_{k=0}^{n-1}\sum_{j=k-m+1}^{k+m-1}\left(E|\bar{h}_k(X_k^m)|^3 E|\bar{h}_j(X_j^m)|^3\right)^{1/3}. \tag{3.14}$$

Now we apply (3.7) and (3.8) to (3.6), and (3.13), (3.14) to (3.9), next we apply the results to (3.4), then we get

$$\int_{-\infty}^{\infty} I(\tau)d\tau = \sqrt{2\pi}(1+o(1))\left(\sum_{k=0}^{n-1}\sum_{j=k-m+1}^{k+m-1}\text{cov}\left(h_k(X_k^m),h_j(X_j^m)\right)-nB_n^2\right)$$

$$+\frac{cm}{n}\sum_{k=0}^{n-1}\sum_{j=k-m+1}^{k+m-1}\left(E|\bar{h}_k(X_k^m)|^3 E|\bar{h}_j(X_j^m)|^3\right)^{1/3}+\frac{cm^{3/2}}{n^{3/2}}\left(\sum_{k=0}^{n-1}(E|\bar{h}_k(X_k^m)|^3)^{1/3}\right)^2. \tag{3.15}$$

Proposition 2.2 follows by (3.15), (3.4), (3.1), definition of $\sigma_{n,m}^2$ (2.1), and Proposition 2.1.

**Proof of Theorems 2.1 and 2.2.** Although Theorem 2.2 is a direct consequence of Theorem 2.1, nevertheless, in cours of the proof of Theorem 2.1 below we will also observe the proof of Theorem 2.2.

In what follows, for simplicity of notation, we will assume $Eh_k(X_k^m)=0, k=0,1,...,n-1$. Let $\{N=N(n)\}$ be a sequence of integers such that $N\leq n-4$ and $N/n\to a\in(0,1)$.

Write

$$R_{n,m}=\sum_{k=0}^{N-m}h_k(nS_k^m)+\sum_{k=N-m+1}^{n-1}h_k(nS_k^m)=:G_{N-m,m}^S+\bar{G}_{N-m,m}^S.$$

Let's consider r.v. $G_{N-m,m}^S$. Set $\tilde{g}_k(x_k^m)=\tilde{h}_k(x_k^m)-\tilde{B}_n(x_k-1)$, where $\tilde{h}_k=h_k$ $k=0,1,...,N-m$, else $\tilde{h}_k=0$, and

$$\tilde{B}_n=\frac{1}{N}\sum_{k=0}^{N-m}\text{cov}\left(h_k(X_k^m),|X_k^m|\right)=B_{N-m,m}(1+o(1)), \tag{3.16}$$

where $B_{N-m,m}$ is defined according to (2.1). Since $S_0+...+S_{n-1}=1$ we have $G_{N-m,m}^S=\sum_{k=0}^{n-1}\tilde{g}_k(nS_k^m)$.

Set $G_{N-m,m}^X=\sum_{k=0}^{n-1}\tilde{g}_k(nS_k^m)$. Evidently $E\tilde{g}_k(X_k^m)=0$, $k=0,1,...,n-1$. We have

$$E\left(G_{N-m,m}^X\right)^2=\sum_{k=0}^{n-1}\sum_{j=0}^{n-1}E\tilde{g}_k(X_k^m)\tilde{g}_j(X_j^m)=\sum_{k=0}^{N-m}\sum_{j=0}^{N-m}Eh_k(X_k^m)h_j(X_j^m)$$

$$-2\tilde{B}_n\sum_{k=0}^{N-m}E\left(h_k(X_k^m)\sum_{j=0}^{n-1}(X_j-1)\right)+\tilde{B}_n^2\sum_{k=0}^{n-1}\sum_{j=0}^{n-1}E(X_k-1)(X_j-1)$$

$$= \sum_{k=0}^{N-m}\sum_{j=k-m+1}^{k+m-1} \text{cov}\left[h_k(X_k^m)h_j(X_j^m)\right] - (2N-n)\tilde{B}_n^2 =: \tilde{\sigma}_{N-m,m}^2 ,$$

Next,

$$\text{Cov}\left(G_{N-m,m}^X, \sum_{l=0}^{N-1}X_l\right) = \sum_{k=0}^{N-m}Eh_k(X_k^m)\sum_{l=0}^{N-1}(X_l-1) - \tilde{B}_n\sum_{l=0}^{N-1}\sum_{l=0}^{n-1}E(X_l-1)(X_l-1)$$

$$=\sum_{k=0}^{N-m}\text{cov}\left(h_k(X_k^m),\left|X_k^m\right|\right) - N\tilde{B}_n = 0 .$$

Thus

$$EG_{N-m,m}^X = 0 , \quad \text{Var}G_{N-m,m}^X = \tilde{\sigma}_{N-m,m}^2 \quad \text{and} \quad \text{Cov}\left(G_{N-m,m}^X, \sum_{l=0}^{N-1}X_l\right) = 0. \qquad (3.17)$$

One can observe that $(2N-n)\tilde{B}_n^2 \sim (n-m)B_{n-m,m}^2$ as $n \to \infty$ and $a \nearrow 1$, and hence by condition (i)

$$\tilde{\sigma}_{N-m,m}^2 \sim \sigma_{n,m}^2, \text{ as } n \to \infty \text{ and } a \nearrow 1 . \qquad (3.18)$$

Set

$$\xi_{N,m}(t,\tau) := t\sigma_{N-m,m}^{-1}G_{N-m,m}^X + \tau n^{-1/2}\sum_{k=0}^{N-1}(X_k-1),$$

$$\psi_{N,m}(t,\tau) := E\exp\{i\xi_{N,m}(t,\tau)\}, \quad \phi_N(\tau) = E\exp\left\{i\tau n^{-1/2}\sum_{k=N}^{n-1}(X_k-1)\right\}.$$

By above Lemma we have

$$\varphi_{N-m}(t) := E\exp\{it\tilde{\sigma}_{N-m,m}^{-1}G_{N-m,m}^S\} = \mu_n^{-1}\int_{-\infty}^{\infty}\psi_{N,m}(t,\tau)\phi_N(\tau)d\tau . \qquad (3.19)$$

Due to (3.17)

$$E\xi_{N,m}(t,\tau) = 0 \quad \text{and} \quad \text{Var}\xi_{N,m}(t,\tau) = t^2 + \tau^2 N/n. \qquad (3.20)$$

Further, $\xi_N(t,\tau)$ is the sum of $m$-dependent r.v.s. With agreement (for the moment) that $\tilde{g}_k(x_k^m) = 0$ for $k = N-m+1,...,N-1$, we obtain for $r > 2$

$$\beta_{N,r}(t,\tau) := \sum_{k=0}^{N-1}E\left|\frac{t}{\tilde{\sigma}_{N-m,m}}\tilde{g}_k(X_k^m) + \frac{\tau}{\sqrt{n}}(X_k-1))\right|^r$$

$$\leq 2^{r-1}|t|^r \tilde{\sigma}_{N-m,m}^{-r}\sum_{k=0}^{N-m}E\left|\tilde{g}_k(X_k^m)\right|^r + 2^{r-1}c(r)|\tau|^r\frac{N}{n^{r/2}},$$

This together with (3.20) yield

$$m^{r-1}\beta_N(t,\tau)/(\text{Var}\tilde{\xi}_N(t,\tau))^{r/2} = m^{r-1}\beta_N(t,\tau)/(t^2+\tau^2N/n)^{r/2}$$

$$\leq c(r)\left(\frac{m^{r-1}}{\tilde{\sigma}_{N-m,m}^r}\sum_{k=0}^{N-m}E\left|\tilde{g}_k(X_k^m)\right|^r + an^{-(r-2)/2}\right) \to 0$$

since $N/n \to a \in (0,1)$, $m = o(n)$ and condition (2.2). So by Theorem 4.1 of Svante (2021) and (3.22) the r.v. $\xi_N(t,\tau)$ has an asymptotically normal distribution whose mean and variance are 0 and $t^2 + a\tau^2$, respectively. Hence for each $t$ and $\tau$

$$\psi_{N,m}(t,\tau) \to \exp\{-(t^2 + a\tau^2)/2\}. \qquad (3.21)$$

By CLT and the fact that $|\phi_N(\tau)| = (1 + \tau^2/n)^{-(n-N)/2}$ one can observe, respectively, that

$$\phi_N(\tau) \to \exp\{-(1-a)\tau^2/2\} \text{ and } \int_{-\infty}^{\infty} |\phi_N(\tau)| d\tau \to \int_{-\infty}^{\infty} \exp\{-(1-a)\tau^2/2\} d\tau. \qquad (3.22)$$

These facts together with (3.21) by Lebesgue dominated convergence theorem yields

$$\varphi_{N-m}(t) \to \frac{1}{\sqrt{2\pi}} \int_{-\infty}^{\infty} \exp\{-(t^2 + a\tau^2)/2\} \exp\{-(1-a)\tau^2)/2\} d\tau = \exp\{-t^2/2\} \text{ as } n \to \infty. \qquad (3.23)$$

By the way, by letting $a \nearrow 1$ in the both sides of (3.23) we conclude the proof of Theorem 2.2. By (3.18) we obtain

$$\lim_{n \to \infty} E \exp\{it\sigma_{n,m}^{-1} G_{n-m,m}^S\} = \lim_{a \nearrow 1} \lim_{n \to \infty} \varphi_{an,m}(t\tilde{\sigma}_{N-m,m}/\sigma_{n,m}) = \exp\{-t^2/2\}. \qquad (3.24)$$

Let's now consider r.v. $\bar{G}_{N-m,m}^S$. Set $\bar{g}_k(x_k^m) = \bar{h}_k(x_k^m) - \bar{B}_n(x_k - 1)$, where $\bar{h}_k = h_k$ $k = N-m,\ldots,n-1$, else $\bar{h}_k = 0$, and $\bar{B}_n = B_{n,m} - (1 - mn^{-1})B_{n-m,m} = o(1)$ by condition (i). We have

$$\bar{G}_{N-m,m}^S = \sum_{k=N-m+1}^{n-1} \bar{g}_k(nS_k^m), \; \bar{G}_{N-m,m}^X = \sum_{k=N-m+1}^{n-1} \bar{g}_k(X_k^m), \; \text{Var}\bar{G}_{N-m,m}^X = \bar{\sigma}_{N-m,m}^2,$$

where $\bar{\sigma}_{N-m,m}^2$ should be rewritten from formula for the $\sigma_{n,m}^2$, see (2.1), by replacing $\bar{h}_k$ instead of $h_k$. Therefore, by condition (i) we obtain

$$\sigma_{n,m}^{-2} \bar{\sigma}_{N-m,m}^2 \to 0 \text{ as } n \to \infty \text{ and } a \nearrow 1. \qquad (3.25)$$

Next,

$$\text{cov}\left(\frac{1}{\sigma_{n,m}} G_{N-m,m}^X, \frac{1}{\sqrt{n}} \sum_{j=N-m+1}^{n-1} (X_j - 1)\right)$$

$$= \frac{\bar{\sigma}_{N-m,m}}{\sigma_{n,m}} \sqrt{1 - \frac{N}{n} + \frac{m-2}{n}} \text{corr}\left(G_{N-m,m}^X, \sum_{j=N-m+1}^{n-1} (X_j - 1)\right) \to 0, \text{as } n \to \infty \text{ and } a \nearrow 1. \qquad (3.26)$$

Also

$$\text{Var} \frac{1}{\sqrt{n}} \sum_{j=N-m+1}^{n-1} (X_k - 1) = (n - N + m - 2)/n \to 1 - a, \text{ as } n \to \infty. \qquad (3.27)$$

Due to Lemma we have

$$\bar{\varphi}_{N-m}(t) := E \exp\{it\sigma_{n,m}^{-1} \bar{G}_{N-m,m}^S\} = \mu_n^{-1} \int_{-\infty}^{\infty} \bar{\psi}_{N,m}(t,\tau) \bar{\phi}_N(\tau) d\tau, \qquad (3.28)$$

where $\eta_{N,m}(t,\tau) := t\sigma_{n,m}^{-1}\bar{G}_{N-m,m}^{X} + \tau n^{-1/2}\sum_{k=N-m+1}^{n-1}(X_k-1)$,

$$\bar{\psi}_{N,m}(t,\tau) := E\exp\{i\eta_{N,m}(t,\tau)\}, \quad \bar{\phi}_N(\tau) = E\exp\left\{i\tau n^{-1/2}\sum_{k=0}^{N-m}(X_k-1)\right\}.$$

Alike to (3.22) we conclude that

$$\bar{\phi}_N(\tau) \to \exp\{-a\tau^2/2\} \text{ and } \int_{-\infty}^{\infty}|\bar{\phi}_N(\tau)|d\tau \to \int_{-\infty}^{\infty}\exp\{-a\tau^2/2\}d\tau. \qquad (3.29)$$

Since (3.25), (3.26) and (3.27) the reasons similar to those carried out to obtain (3.23) allows to get

$$\bar{\psi}_{N,m}(t,\tau) \to \exp\{-d_1(a)t^2/2 + 2t\tau d_2(a) - (1-a)\tau^2/2\},$$

where $d_1(a) \to 0$ and $d_2(a) \to 0$ as $a \nearrow 1$. This together with (3.1), (3.28), (3.29) imply

$\bar{\varphi}_{N,m}(t) \to 1$, as $n \to \infty$ and $a \nearrow 1$. Hence

$$\lim_{n\to\infty}E\exp\{it\sigma_{n,m}^{-1}\bar{G}_{n-m,m}^{S}\} = \lim_{a\nearrow 1}\lim_{n\to\infty}\bar{\varphi}_{N,m}(t) = 1.$$

This and (3.24) yield by Lemma 5 of Le Cam (1958) $E\exp\{it\sigma_{n,m}^{-1}G_{n,m}^{S}\} \to \exp\{-t^2/2\}$. Theorem 2.1 follows.

**Proofs of Corollaries 2.1, 2.2, 2.3** and **2.4** follow in immediate manner from Propositions 2.1, 2.2 and Theorem 2.1, respectively.

**Proof of Corollaries 2.5**. We have $E|X_0^m|^2 = m(m+1)$ and letting $h(x) = x^2$ we obtain

$\sigma_m^2 = 2m(m+1)(2m+1)/3 = 4m^3/3(1+O(m^{-1}))$ as $m \to \infty$. On the other hand

$$\sum_{k=0}^{n-1}(n|S_k^m|)^2 - nm(m+1) = \sum_{k=0}^{n-1}(n|S_k^m|-m)^2 - nm,$$

since $S_{n+j} = S_j$, $j \geq 0$, and hence $|S_0^m|+...+|S_{n-1}^m| = m$. Hence we can take $h(u) = (u-m)^2$ and

$g(|X_0^m|) = (|X_0^m|-m)^2 - m(2X_0-1)$. Then $E|g(|X_0^m|)|^3 \leq cm^3$. Part (i) follows.

Let $h(u) = \ln x$. By Abramowitz and Stegun 1972, p. 261, and Gradshteyn and Ryzhik (2000, p. 577):

$$\psi(m+1) = \psi(m) + m^{-1}, \quad \psi(m) = \ln m - \frac{1}{2m} - \frac{1}{12m^2} + O(m^{-3}), \quad \zeta(2,m) = \frac{1}{m} + \frac{1}{2m^2} + O(m^{-3}). \quad (3.30)$$

We have $E\ln|X_k^m| = \psi(m)$ and $\text{cov}(\ln|X_0^m|,|X_0^m|) = m(\psi(m+1)-\psi(m)) = 1$. Also, see Holst (1979), Kuo and Jammalamadaka (1981), $\sigma_m^2 = (2m^2-2m+1)\zeta(2,m) - 2m+1$. Next,

$g(|X_0^m|) = \ln|X_0^m| - \psi(m) - (X_0-1)$. Using (3.30) after some lengthy calculations we find

$\sigma_m^2 = (2m^2)^{-1}(1+O(m^{-1}))$ and $E[g(|X_0^m|)]^4 = O(m^{-4})$ as $m \to \infty$. Eq (2.6) follows.

Let now $h(x) = x \ln x$, with convenience $0 \ln 0 = 0$. We have $E|X_0^m| \ln |X_0^m| = m\psi(m+1)$ and $\text{cov}(|X_0^m| \ln |X_0^m|, |X_0^m|) = m(1 + \psi(m+1))$. By Misra (2001, p.1452)

$\sigma_m^2 = (2(m(m+1))^2 \zeta(2, m+2) - m(m+1)(2m-1))/4$. Using (3.30) after some calculations we get

$\sigma_m^2 = (m+5)/4 + O(m^{-1})$. Further $g(|X_0^m|) = |X_0^m| \ln |X_0^m| - m\psi(m+1) - (X_0 - 1)m(1 + \psi(m+1))$. By formulas (4.353) and (4.358) of Gradshteyn and Ryzhik (2000, pp. 574, 576,577 ) and asymptotic expansion formula for $\zeta(k, m)$, cf.Abramowitz and Stegun (1972, p. 261) after some too lengthy calculations we obtain $E[g(|X_0^m|)]^4 = O(m^2)$. Eq (2.7) concludes by Corollary 2.4.

**Conflict of interests**

The author has no relevant financial or non-financial interests to disclose.